\documentclass[12pt, reqno]{conm-p-l}
\usepackage[dvips]{graphicx}

\numberwithin{equation}{section} \numberwithin{figure}{section}
\numberwithin{table}{section}

\DeclareMathOperator{\fillrad}{{\rm FillRad}} \DeclareMathOperator{\fillvol}{{\rm
FillVol}} \DeclareMathOperator{\length}{{\rm length}}
\DeclareMathOperator{\vol}{{\rm vol}}

\DeclareMathOperator{\area}{{\rm area}} 

\DeclareMathOperator{\rk}{{\rm rank}} \DeclareMathOperator{\PD}{{\rm PD}}

\DeclareMathOperator{\stsys}{{\rm stsys}} \DeclareMathOperator{\sys}{{\rm sys}}
\DeclareMathOperator{\confsys}{{\rm confsys}}

\DeclareMathOperator{\diam}{{\rm diam}} \DeclareMathOperator{\XX}{{\overline{X}}}
\DeclareMathOperator{\fx}{{[\overline{F_{\it \! X}}]}}

\DeclareMathOperator{\pisys}{{\rm sys}\pi}

\DeclareMathOperator{\AJ} {{\mathcal A}}

\newcommand{\C}{\mathbb C} 
\newcommand{\R}{\mathbb R} \newcommand{\T}{\mathbb T}
\newcommand{\Z}{\mathbb Z}

\newtheorem{theorem}{Theorem}[section]
\newtheorem{proposition}[theorem]{Proposition}
\newtheorem{lemma}[theorem]{Lemma}
\newtheorem{corollary}[theorem]{Corollary}
\newtheorem{conjecture}[theorem]{Conjecture}

\theoremstyle{definition}
\newtheorem{definition}[theorem]{Definition}
\newtheorem{example}[theorem]{Example}
\newtheorem{remark}[theorem]{Remark}

\begin{document}

\author[M.~Katz]{Mikhail G. Katz$^{*}$} \address{Department of
Mathematics \\ Bar-Ilan University\\ Ramat-Gan 52900 \\Israel}
\email{katzmik@math.biu.ac.il} \thanks{$^{*}$Supported by the Israel Science
Foundation (grants no.\ 620/00-10.0 and 84/03)}

\author[C.~Lescop]{Christine Lescop$^{\dagger}$} \address{Institut Fourier\\ BP 74, Universit\'e de Grenoble I\\ 38402 Saint Martin
d'H\`eres cedex\\ France} \email{Christine.Lescop@ujf-grenoble.fr}
\thanks{$^{\dagger}$Supported by CNRS (UMR 5582)}

\title[Filling Area Conjecture and the Fiber Class] {Filling Area
Conjecture, Optimal Systolic Inequalities, and the Fiber Class in Abelian Covers}

\dedicatory{Dedicated to the memory of Robert Brooks, a colleague and friend.}


\copyrightinfo{2005}{M. Katz, C. Lescop}


\begin{abstract}
We investigate the filling area conjecture, optimal systolic
inequalities, and the related problem of the nonvanishing of certain
linking numbers in 3-manifolds.
\end{abstract}

\maketitle

\tableofcontents

\section{Filling radius and systole}

\begin{figure}[h]
\includegraphics [width=8cm] {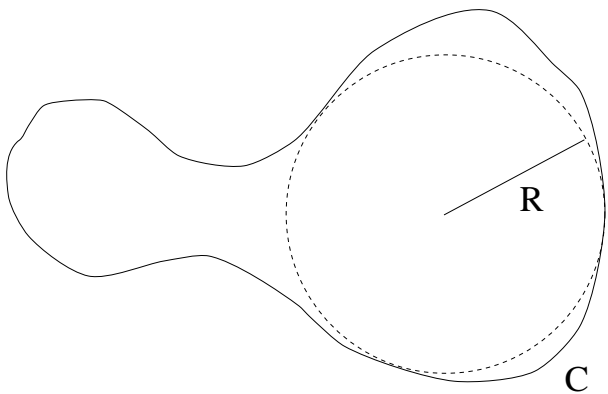} \caption{Largest inscribed circle has
radius $R$} \label{11b}
\end{figure}

The filling radius of a simple loop $C$ in the plane is defined as the largest
radius, $R>0$, of a circle that fits inside $C$, see Figure~\ref{11b}:
\[
\fillrad(C\subset \R^2) = R.
\]

\begin{figure}[h]
\includegraphics[width=10cm]{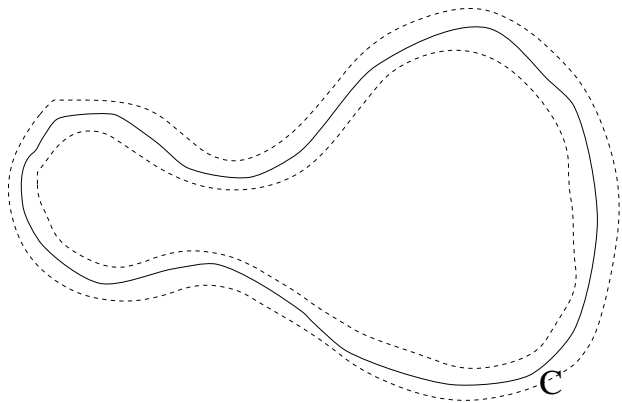} \caption{Neighborhoods
of loop $C$} \label{12}
\end{figure}

There is a kind of a dual point of view that allows one to generalize this notion in
an extremely fruitful way, as shown in \cite{Gr1}. Namely, we consider the
$\epsilon$-neighborhoods of the loop $C$, denoted $U_\epsilon C \subset \R^2$, see
Figure \ref{12}.

As $\epsilon>0$ increases, the $\epsilon$-neighborhood $U_\epsilon C$ swallows up
more and more of the interior of the loop.  The ``last'' point to be swallowed up is
precisely the center of a largest inscribed circle.  Therefore we can reformulate
the above definition by setting
\[
\fillrad(C\subset \R^2) = \inf \left\{ \epsilon > 0 \left| \mbox{ loop $C$ contracts
to a point in $U_\epsilon C$} \right. \right\}.
\]

Given a compact manifold $X$ imbedded in, say, Euclidean space $E$, we could
define the filling radius {\em relative\/} to the imbedding, by minimizing the
size of the neighborhood $U_\epsilon X \subset E$ in which $X$ could be
homotoped to something smaller dimensional, e.g., to a lower dimensional
polyhedron. Technically it is more convenient to work with a homological
definition.

Denote by $A$ the coefficient ring $\Z$ or $\Z_2$, depending on whether or not $X$
is oriented.  Then the fundamental class, denoted $[X]$, of a compact
$n$-dimensional manifold $X$, is a generator of $H_n(X;A)=A$, and we set
\begin{equation}
\label{11e} \fillrad(X\subset E) = \inf \left\{ \epsilon > 0 \left|
\;\iota_\epsilon^{\phantom{I}}([X])=0\in H_n(U_\epsilon X) \right. \right\},
\end{equation}
where $\iota_\epsilon$ is the inclusion homomorphism.

To define an ``absolute'' filling radius in a situation where $X$ is
equipped with a Riemannian metric~$g$, M. Gromov proceeds as follows,
cf.  \eqref{fillrad}. One exploits an imbedding due to C. Kuratowski
\cite{Kur} (while the spelling of the author's first name used in
\cite{Kur} starts with a~``C'', the correct spelling seems to be
Kazimierz).  One imbeds $X$ in the Banach space $L^\infty(X)$ of
bounded Borel functions on $X$, equipped with the sup
norm~$\|\;\|$. Namely, we map a point $x\in X$ to the function~$f_x\in
L^\infty(X)$ defined by the formula $f_x(y) = d(x,y)$ for all $y\in
X$, where $d$ is the distance function defined by the metric. By the
triangle inequality we have
\begin{equation}
\label{11c} d(x,y) = \| f_x - f_y \| ,
\end{equation}
and therefore the imbedding is strongly isometric, in the precise sense of
formula~\eqref{11c} which says that internal distance and ambient distance coincide.
Such a strongly isometric imbedding is impossible if the ambient space is a Hilbert
space, even when $X$ is the Riemannian circle (the distance between opposite points
must be $\pi$, not 2!).  We then set~$E= L^\infty(X)$ in formula~\eqref{11e} and
define
\begin{equation}
\label{fillrad} \fillrad(X)=\fillrad \left( X\subset L^{\infty\phantom{I}}
_{\phantom{I}}\!\!\!(X) \right) .
\end{equation}
Exact calculation of the filling radius is discussed in Section~\ref {nine}.  The
related invariant called filling volume is defined in Section~\ref{four}.  For the
relation between filling radius and systole, see \eqref{91b}.  The defining text for
this material is \cite{Gr3}, with more details in \cite{Gr1,Gr2}.  See also the
recent survey \cite{CK}.

The systole, $\pisys_1(g)$, of a compact non simply connected Riemannian
manifold~$(X,g)$ is the least length of a noncontractible loop $\gamma\subset X$:
\begin{equation}
\label{ps1} \pisys_1(g)=\min_{[\gamma]\not= 0\in \pi_1(X)}\length(\gamma).
\end{equation}

This notion of systole is apparently unrelated to the systolic arrays
of~\cite{Ku}.  We will be concerned with comparing this Riemannian
invariant to the total volume of the metric, as in Loewner's
inequality~\eqref{(1.1)}.

In Section~\ref{one}, we recall the classical Loewner inequality, as well as
Pu's inequality.  In Section~\ref{four}, we define an invariant called the
filling volume, and state Gromov's filling volume (e.g., area) conjecture.  The
notion of the successive minima of a lattice in Banach space is recalled in
Section~\ref{five}.  The results on the precise calculation of filling
invariants (both the filling radius and the filling area) are reviewed in
Section~\ref {nine}.  Stable and conformal systoles are defined in
Section~\ref{six}.

In Section~\ref{seven}, we recall the definition of the Abel-Jacobi
map to the torus, and present M. Gromov's optimal stable systolic
inequality for $n$-tori.  Higher dimensional optimal generalisations
of Loewner's inequality are discussed in Section~\ref{bank} and
Section~\ref{ik}.  Thus, in Section~\ref{bank}, we recall J. Hebda's
inequality combining systoles of dimension and codimension 1 in the
case of unit Betti number, and present a generalisation to arbitrary
Betti number.

Optimal inequalities combining varieties of 1-systoles appear in
Section~\ref{ik}.  Here a necessary topological condition is the nonvanishing
of the homology class of the lift to the maximal free abelian cover, of the
typical fiber of the Abel-Jacobi map.  The relation of this condition to the
generalized Casson invariant $\lambda$ of \cite{Les} is examined in
Section~\ref{thirteen}.  The example of nilmanifolds is discussed in
Section~\ref{eleven}, and the general case treated in Section~\ref{marin}.

\section{Inequalities of C. Loewner and P. Pu}
\label{one}

The Hermite constant, denoted $\gamma_n$, can be defined as the
optimal constant in the inequality
\begin{equation}
\label{11} \pisys_1(\T^n)^2 \leq \gamma_n \vol(\T^n)^{2/n},
\end{equation}
over the class of all {\em flat\/} tori $\T^n$, cf. Section~\ref {five}.  Here
$\gamma_n$ is asymptotically linear in $n$, cf. \cite[pp.~334, 337]{LLS}.  The
precise value is known for small $n$, e.g., $\gamma_2= \frac{2}{\sqrt{3}}$ (see
Lemma~\ref{217}), $\gamma_3 = 2^{\frac{1}{3}}, \ldots$.  It can be shown via
the filling invariants that an inequality of type~\eqref{11} remains valid in
the class of {\em all\/} metrics, but with a nonsharp constant on the order of
$n^{4n}$ \cite{Gr1}.  M. Gromov's proof uses inequality~\eqref{91b} as a
starting point.

Around 1949, Charles Loewner proved the first systolic inequality, cf.
\cite{Pu}. He showed that every Riemannian metric $g$ on the torus~$\T^2$
satisfies the inequality
\begin{equation}
\label{(1.1)} \pisys_1(g)^2\le \gamma_2 \area(g),
\end{equation}
while a metric satisfying the boundary case of equality in \eqref{(1.1)} is
necessarily flat, and is homothetic to the quotient of $\C$ by the lattice spanned
by the cube roots of unity.

The following estimate is found in \cite[Corollary~5.2.B]{Gr1}. Namely, every
aspherical compact surface $(\Sigma,g)$ admits a metric ball
$B=B_p\left(\tfrac{1}{2}\pisys_1(g)\right) \subset \Sigma$ of radius
$\tfrac{1}{2}\pisys_1(g)$, which satisfies

\begin{equation}
\label{tqi} \pisys_1(g)^2 \leq \frac{4}{3}\area(B).
\end{equation}

P. Pu's inequality \cite{Pu} admits an immediate generalisation via Gromov's
inequality \eqref{tqi}.  Namely, every surface $(X,g)$ which is not a 2-sphere
satisfies the inequality
\begin{equation}
\label{pu3} \pisys_1(g)^2\leq \frac{\pi}{2}\area(g),
\end{equation}
where the boundary case of equality in \eqref{pu3} is attained precisely when, on
the one hand, the surface $X$ is a real projective plane, and on the other, the
metric $g$ is of constant Gaussian curvature.

\section{Filling area conjecture}
\label{four}

Consider a compact manifold $N$ of dimension $n \geq 1$ with a distance
function~$d$.  Here $d$ could be induced by a Riemannian metric~$g$, but could also
be more general.  The notion of the Filling Volume, $\fillvol(N^n,d)$, of $N$ was
introduced in~\cite{Gr1}, where it is shown that when $n \geq 2$,
\begin{equation}
\label{fillvoldef} \fillvol(N^n,d)=\inf_g{\vol_{n+1}(X^{n+1},g)}
\end{equation}
where $X$ is any fixed manifold such that $\partial X=N$.  Here one can even
take a cylinder $X=N\times [0,\infty)$.  The infimum is taken over all complete
Riemannian metrics $g$ on $X$ for which the boundary distance function is $\geq
d$, i.e., the length of the shortest path in $X$ between boundary points $p$
and $q$ is $\geq d(p,q)$.  In the case $n=1$, the topology of the 2-dimensional
filling can affect the infimum, as is shown by example in
\cite[Counterexamples~2.2.B]{Gr1}.

The precise value of the filling volume is not known for any Riemannian metric.
However, the following values for the canonical metrics on the spheres (of sectional
curvature $+1$) is conjectured in \cite{Gr1}, immediately after Proposition 2.2.A.
\begin{conjecture}
\label{fillsphere} $\fillvol(S^n,can)=\frac 1 2 \omega_{n+1}$, where $\omega_{n+1}$
denotes the volume of the unit $(n+1)$-sphere.
\end{conjecture}
This conjecture is still open in all dimensions.  The case $n=1$ can be broken
into separate problems depending on the genus $s$ of the filling, cf.
\eqref{42}.

\section{Lattices and successive minima}
\label{five}

Let $b>0$.  By a {\em lattice\/} $L$ in Euclidean space $\R^b$, we mean a
discrete subgroup isomorphic to $\Z^b$, i.e., the integer span of a linearly
independent set of $b$ vectors.

More generally, let $B$ be a finite-dimensional Banach space, i.e., a vector
space together with a norm $\|\;\|$.  Let $L\subset (B,\|\;\|)$ be a lattice of
maximal rank $\rk(L)=\dim(B)$.

\begin{definition}
\label{mindef} For each $k=1,2,\dots, \rk(L)$, define the {\em $k$-th successive
minimum} of the lattice $L$ by
\begin{equation}
\label{success} \lambda_k(L,\|\;\|) = \inf\bigg\{\lambda\in\R \ \left| \,
\begin{array}{l}%
\exists \text{ lin. indep. } v_1, \ldots, v_k \in L \\ \text{\ with\ } \|v_i\|\leq
\lambda
\end{array}
\right. \bigg\}.
\end{equation}
\end{definition}

Thus, the first successive minimum, $\lambda_1(L,\|\;\|)$ is the least length of a
nonzero vector in $L$.  If $\rk(L)\geq 2$, then given a pair of vectors $S=\{v, w\}$
in $L$, define the length $|S|$ of $S$ by setting $|S|=\max(\|v\|, \|w\|)$.  Then
the second successive minimum, $\lambda_2(L,\|\;\|)$ is the least length of a pair
of non-proportional vectors in $L$: $\lambda_2(L)= \inf_S |S|$, where $S$ runs over
all linearly independent pairs $\{v, w\}\subset L$.

\begin{example}
The standard lattice $\Z^b\subset \R^b$ satisfies
\[
\lambda_1(\Z^b)= \ldots=\lambda_b(\Z^b)=1.
\]
Note that the torus
$\T^b= \R^b/\Z^b$ satisfies $\vol(\T^b)= 1$ as it has the unit cube as
a fundamental domain.
\end{example}

Here the volume of the torus $\R^b/L$ (also called the covolume of the
lattice~$L$) is by definition the volume of a fundamental domain for $L$, e.g.,
a parallelepiped spanned by a~$\Z$-basis for~$L$.  The volume can be calculated
as the square root of the determinant of the Gram matrix of such a basis. Given
a finite set $S= \{v_i\}$ in $\R^b$, we define its {\em Gram matrix\/} as the
matrix of inner products $Gram(S)= (\langle v_i, v_j \rangle)$.  Then the
parallelepiped $P$ spanned by the vectors $\{v_i\}$ satisfies $\vol(P)= {\det (
Gram(S))}^{1/2}.$

\begin{example}
Consider the lattice $L_\zeta \subset \R^2=\C$ spanned by $1\in \C$ and the sixth
root of unity $\zeta=e^{\frac{2\pi i}{6}}\in \C$.  Then $L_\zeta = \Z \zeta + \Z 1$
satisfies $\lambda_1 (L_\zeta) = \lambda_2(L_\zeta)=1$.  Meanwhile, the torus
$\T^2=\R^2/L_\zeta$ satisfies $\area(\T^2)=\frac{\sqrt{3}}{2}$.
\end{example}

\begin{definition}
\label{53} The Hermite constant $\gamma_b$ is defined in one of the following two
equivalent ways:
\begin{enumerate}
\item[(a)]
$\gamma_b$ is the {\em square\/} of the maximal first successive minimum, among all
lattices of unit covolume;
\item[(b)]
$\gamma_b$ is defined by the formula
\begin{equation}
\label{???} \sqrt{\gamma_b} = \sup\left\{ \left. \frac{\lambda_1(L)}{\vol(\R^b/L)
^{\frac{1}{b}}} \right| L \subseteq (\R ^b, \|\;\|) \right\},
\end{equation}
where the supremum is extended over all lattices $L$ in $\R^b$ with a Euclidean norm
$\|\;\|$.
\end{enumerate}
\end{definition}

A lattice realizing the supremum is called {\em critical}.  Such a lattice may be
thought of as realizing a densest packing in $\R^b$, when we place balls of radius
$\frac{1}{2}\lambda_1(L)$ at the points of $L$.

\begin{lemma}
\label{217} Let $b=2$.  Then we have the following value for the Hermite constant:
$\gamma_2=\frac{2}{\sqrt{3}}=1.1547...$.  The corresponding optimal lattice is
homothetic to the $\Z$-span of cube roots of unity in $\C$.
\end{lemma}

\begin{proof}
Consider a lattice $L\subset \C=\R^2$.  Clearly, multiplying $L$ by nonzero
complex numbers does not change the value of the quotient
$\frac{\lambda_1(L)^2}{\area(\C/L)}.$ Choose a ``shortest'' vector $z \in L$,
i.e., we have $|z|= \lambda_1 (L)$.  By replacing $L$ by the lattice $z^{-1}L$,
we may assume that the complex number $+1\in \C$ is a shortest element in the
lattice.  We will denote the new lattice by the same letter $L$, so that now
$\lambda_1(L)=1$.  Now complete~$+1$ to a $\Z$-basis $\{\tau, 1\}$, so that
\[
L= \Z\tau + \Z 1.
\]
Thus $|\tau|\geq \lambda_1(L)=1$.  Consider the real part $\Re(\tau)$. We can
adjust the basis by adding or subtracting a suitable integer to $\tau$, so as
to satisfy the condition $-\frac{1}{2} \leq \Re(\tau)\leq \frac{1}{2}.$ Thus,
the second basis vector $\tau$ may be assumed to lie in the closure of the
standard fundamental domain, cf. \cite[p.~78]{Se}
\[
D=\left\{ z\in\C \left|\; |z| > 1, \; | \Re(z) | < \frac{1}{2}, \Im(z)>0 \right.
\right\}
\]
(for the action of $PSL(2,\Z)$ in the upper half-plane of $\C$).  Then the imaginary
part satisfies $\Im(\tau)\geq \frac{\sqrt{3}}{2}$, with equality possible precisely
for~$\tau=e^{i\frac{\pi}{3}} \hbox{\ or\ } \tau=e^{i\frac {2\pi}{3}}.$ The proof is
concluded by calculating the area of the parallelogram in $\C$ spanned by $\tau$ and
$+1$ as follows: $\frac{\lambda_1(L)^2} {\area(\C/L)} = \frac{1}{\Im(\tau)} \leq
\frac{2}{\sqrt{3}}.$
\end{proof}

\begin{example}
In dimensions $b\geq 3$, the Hermite constants are harder to compute, but
explicit values (as well as the associated critical lattices) are known for
small dimensions, e.g., $\gamma_3=2^{\frac{1}{3}} =1.2599...$, while in
dimension 4, one has~$\gamma_4=\sqrt{2}=1.4142...$ \cite{CS}.
\end{example}

The lattice $L^*$ dual to $L$ may be described as follows.  If $L$ is the $\Z$-span
of vectors $(x_i)$, then $L^*$ is the $\Z$-span of a dual basis $(y_j)$ satisfying
\begin{equation}
\label{4.3} \langle x_i, y_j \rangle = \delta_{ij}.
\end{equation}

\begin{definition}
The Berg\'e-Martinet constant, denoted $\gamma_b'$, is defined as follows:
\begin{equation}
\label{22'} \gamma'_b = \sup\left\{ \lambda_1(L) \lambda_1(L^*)\left| L \subseteq \R
^b \right. \right\},
\end{equation}
where the supremum is extended over all lattices $L$ in $\R^b$ with its Euclidean
norm.  A lattice attaining the supremum is called {\em dual-critical}.
\end{definition}

Like the Hermite constant $\gamma_b$, the Berg\'e-Martinet constant $\gamma'_b$ is
asymptotically linear in $b$.  Its value is known in dimensions up to~4.

\begin{example}
\label{22bis} In dimension 3, the value of the Berg\'e-Martinet constant,
\linebreak $\gamma _3'= \sqrt{\frac{3}{2}}=1.2247\ldots$, is slightly below the
Hermite constant $\gamma_3= 2^{\frac{1}{3}}=1.2599\ldots$.  It is attained by
the face-centered cubic lattice, which is not isodual
\cite[p.~31]{MH},\linebreak \cite[Proposition 2.13(iii)]{BM}, \cite{CS}.
\end{example}

\section{Precise calculation of filling invariants}
\label{nine}

Unlike the filling area and volume, the filling radius lends itself somewhat easier
to precise calculation.  The following result was proved in~\cite{Ka0}.

\begin{theorem}
\label{91} Let $X$ be a length space with distance function $d$.  Let $Y\subset X$
be a subset.  Let $R>0$.  Assume that the following two conditions are satisfied:
\begin{enumerate}
\item[(a)]
$\diam(Y) \leq 2R$;
\item[(b)]
$d(x,Y)\leq 2R$ for all $x\in X$.
\end{enumerate}
Then $\fillrad(X) \leq R$.
\end{theorem}

\begin{corollary}
Let $S^1$ be the Riemannian circle.  Then $$\fillrad(S^1)=
\frac{1}{6}\length(S^1).$$ More precisely, the neighborhood $U_\epsilon S^1 \subset
L^\infty$ has the homotopy type of $S^1$ for $\epsilon < \frac{1}{2} d_1$, and the
homotopy type of the sphere~$S^3$ when $\frac{1}{2} d_1 < \epsilon < \frac{1}{2}
d_2$, where $d_i= \frac{i}{2i+1}\length(S^1)$.
\end{corollary}

\begin{proof}
Using the set of vertices of a regular inscribed triangle as our subset $Y\subset
S^1$ in Theorem~\ref{91} proves the upper bound.  The lower bound follows from the
general inequality
\begin{equation}
\label{91b} \pisys_1(X) \leq 6\fillrad(X)
\end{equation}
valid for an arbitrary essential manifold $X$ \cite{Gr1}.  The more precise result
on the homotopy type of the neighborhood was proved in \cite{Katz0} using a kind of
Morse theory on the space of subsets of $S^1$, for which regular odd polygons are
the extrema, as discussed below.
\end{proof}

Let $d_i(X)$ be the $i$-th extremal value of the diameter functional on the space of
finite subsets of $X$, with the convention that $d_0=0$.  Here a subset $Y\subset X$
is extremal if small perturbations of size $\epsilon>0$ can decrease the diameter of
$Y$ at most quadratically in $\epsilon$.  For the Riemannian circle we have a
complete classification of the extrema, namely sets of vertices of odd regular
polygons inscribed in the circle.  Thus, we obtain the following extremal values:
\begin{equation}
d_i(S^1) = \frac{i}{2i+1}\length(S^1).
\end{equation}

For the $n$-sphere, the first extremum is the set of vertices of a regular inscribed
$n+1$-simplex, but the higher ones are difficult to classify even on $S^2$
\cite{Ka89}.

The precise calculation of the filling radius is possible in the following cases
\cite{Ka0, Katz0, Katz1, Katz2}.

\begin{theorem}
The filling radius of the following two-point homogeneous spaces $X$
equals half the first extremal value $d_1(X)$:
\begin{enumerate}
\item[(a)] $X=S^n$, where $d_1(S^n) = \arccos \left( - \frac{1}{n+1}
\right)K^{-1/2}$, where $K$ denotes the (constant) value of sectional
curvature;
\item[(b)]
$X=\R P^n$, where $d_1(\R P^n)= \frac{\pi}{3} K^{-1/2}$;
\item[(c)]
$X=\C P^2$, where $d_1(\C P^2) = \arccos \left( - \frac{1}{3} \right)K^{-1/2}$,
where $K$ denotes the maximal value of sectional curvatures.
\end{enumerate}
Meanwhile, for $\C P^3$ we have a strict inequality
\[
\fillrad(\C P^3) > \frac{1}{2} \; d_1(\C P^3) = \frac{1}{2}\arccos \left( -
\tfrac{1}{3} \right)K^{-1/2}.
\]
\end{theorem}

\begin{remark}
In general, understanding the homotopy type of neighborhoods of $X$ in
$L^\infty$ depends on the study of the extremal values of the diameter
functional, which is an interesting open question.  For example, the
value of $d_2(\C P^n)$ is unknown for $n\geq 2$, but progress was made
in \cite{Katz2}.  Successful calculation of these extrema depends on
analyzing certain semialgebraic sets, defined in terms of the distance
functions among points, in the Cartesian (or symmetric) powers of the
space.
\end{remark}

Recently, L. Liu \cite{Liu} undertook the study of the mapping
properties of the filling radius.

We now turn to the related invariant of the filling area.  We report
on the recent progress on calculating the optimal ratio $\sigma_{2,s}$
of \eqref{42}, in the case of hyperelliptic fillings, in the following
sense.  We introduce the constant
\begin{equation}
\label{42} 
\sigma_{2,s} = \sup_{\Sigma_s} \frac{\length(S^1)^2}{\area(\Sigma_s)},
\end{equation}
where $\Sigma_s$ denotes an orientable filling, of genus $s$, of a
Riemannian circle~$S^1$, meaning that the inclusion of the boundary
$S^1 \subset \Sigma_s$ is strongly isometric,
cf. formu\-la~\eqref{11c}. A related invariant appears in the systolic
inequality~\eqref{ptl}.

For simply connected fillings (i.e., by a disk $D$), we have $\sigma _{2,0}
=2\pi$, by Pu's inequality~\eqref{pu3} applied to the real projective plane
obtained from~$\Sigma_0=D$ by identifying pairs of opposite points of the
boundary circle. The case of genus 1 fillings has recently been solved, as
well~\cite{BCIK1}.  Thus,
\begin{equation}
\label{43} \sigma_{2,s}=2\pi \mbox{ for $s=0$ or $1$}.
\end{equation}

Given a filling by an orientable surface~$X_1$ of genus $s$, one identifies
pairs of antipodal points of the boundary, so to form a nonorientable closed
surface $X_2$. We construct the orientable double cover $X_3\to X_2$ with deck
transformation~$\tau$.  The pair $(X,\tau)$ is thus a real Riemann surface of
genus $2s$, which can be called ovalless since $\tau$ is fixed point free, cf.
\cite{BCIK1}. Assuming that $X_3$ is in addition hyperelliptic \cite{Mi}, one
can prove that it satisfies a bound equivalent to the formula
$\sigma_{2,s}=2\pi$.  The result follows from a combination of two ingredients.
On the one hand, one exploits integral geometric comparison with orbifold
metrics of constant positive curvature on real surfaces of even positive genus,
where the conical singularities correspond to Weierstrass points on the
hyperelliptic surface.  On the other hand, one exploits an analysis of the
combinatorics on unions of closed curves, arising as geodesics of such orbifold
metrics.

\section{Stable and conformal systoles}
\label{six}

In a Riemannian manifold $(X,g)$, we define the volume $\vol_k(\sigma)$ of a
Lipschitz $k$-simplex $\sigma: \Delta^k \rightarrow X$ to be the integral over the
$k$-simplex~$\Delta^k$ of the ``volume form'' of the pullback $\sigma^*(g)$.  The
stable norm $\|h\|$ of an element $h\in H_k(X;\R)$ is the infimum of the
volumes~$\vol_k(c)=\Sigma_i |r_i| \vol_k(\sigma_i)$ over all real Lipschitz cycles
$c=\Sigma_i r_i \sigma_i$ representing $h$.  We define the stable $k$-systole of the
metric $g$ by setting
\begin{equation}
\label{41b} \stsys_k(g) =\lambda_1 \left( H_k^{\phantom{I}}(X; \Z)_\R, \|\;\|
\right) ,
\end{equation}
where $\|\;\|$ is the stable norm in homology associated with the metric~$g$,
while $H_k^{\phantom{I}}(X; \Z)_\R$ denotes the integral lattice in real
homology, and $\lambda_1$ is the first successive minimum, cf.
Definition~\ref{mindef}.

The conformally invariant norm $\|\;\|_{n}$ in $H_1(X^n;\R)$ is by definition dual
to the conformally invariant $L^n$-norm in de Rham cohomology.  The latter norm on
$H^1(X,\R)$ is the quotient norm of the corresponding norm on closed one-forms.
Thus, given $\alpha\in H^1(X,\R)$, we set
\[
\|\alpha\|_{n} = \displaystyle \inf_\omega \left\{ \left. \left ( \int_X
|\omega_x|^n {\it dvol}(x) \right)^{\frac{1}{n}} \right| \omega\in \alpha \right\},
\]
where $\omega$ runs over closed one-forms representing $\alpha$.  The conformal
1-systole of $(X^n,g)$ is the quantity
\begin{equation}
\label{57} \confsys_1(g)=\lambda_1\left(H_1(X;\Z)_\R, \|\;\|_{n}\right),
\end{equation}
satisfying $\stsys_1(g) \leq \confsys_1(g)\vol(g)^{\frac{1}{n}}$ on an $n$-manifold
$(X,g)$.

\section{Gromov's optimal inequality for \boldmath$n$-tori}
\label{seven}

Let $X$ be a smooth compact manifold with positive first Betti number. We define the
Jacobi variety $J_1(X)$ of $X$ by setting
\begin{equation}
\label{26} J_1(X)=H_1(X;\R)/H_1(X;\Z)_\R .
\end{equation}
The homotopy class of the {\em Abel-Jacobi map\/} can be defined as follows.
Consider the abelianisation homomorphism on the fundamental group of $X$, followed
by quotienting by torsion.  Consider a map $K(\pi_1(X),1) \to \T^{b_1(X)}=J_1(X)$
associated with the resulting homomorphism.  The composition $X \to K(\pi_1(X),1)
\to \T^{b_1(X)}$ represents the desired homotopy class.  When $X$ is equipped with a
metric, the Abel-Jacobi map
\begin{equation}
\label{ajm} \AJ_X: X \to J_1(X),
\end{equation}
can be induced by the harmonic one-forms on $X$, originally introduced by
A.~Lichnerowicz~\cite{Li}, cf. \cite[4.21]{Gr3}.  More precisely, let $E$ be
the space of harmonic $1$-forms on $X$, with dual $E^*$ canonically identified
with $H_1(X,\R)$.  By integrating an integral harmonic $1$-form along paths
from a basepoint $x_0\in X$, we obtain a map to $\R/\Z=S^1$. In order to define
a map $X\to J_1(X)$ without choosing a basis for cohomology, we argue as
follows. Let $x$ be a point in the universal cover $\widetilde{X}$ of $X$.
Thus~$x$ is represented by a point of $X$ together with a path $c$ from~$x_0$
to it.  By integrating along the path~$c$, we obtain a linear form, $h\to
\int_c h$, on~$E$.  We thus obtain a map $\widetilde{X}\to E^* = H_1(X,\R)$,
which, furthermore, descends to a map
\begin{equation}
\label{72} \overline{\AJ}_X: \overline{X}\to E^*,\;\; c\mapsto \left( h\mapsto
\int_c h \right),
\end{equation}
where $\overline{X}$ is the universal free abelian cover.  By passing to quotients,
this map descends to the Abel-Jacobi map \eqref{ajm}. The Abel-Jacobi map defined
using harmonic forms is unique up to translations of the Jacobi torus $J_1(X)$.

An optimal higher-dimensional generalisation of Loewner's
inequality~\eqref{(1.1)} is due to M.~Gromov \cite[pp.~259-260]{Gr3} (cf.
\cite[inequality (5.14)]{CK}) based on the techniques of D. Burago and S.
Ivanov~\cite{BI94, BI95}.  We give below a slight generalisation of Gromov's
statement.

\begin{definition}
Given a map $f:X\to Y$ between closed manifolds of the same dimension, we denote by
$\deg(f)$, either the algebraic degree of~$f$ when both manifolds are orientable, or
the absolute degree~\cite{BBS}, otherwise.
\end{definition}

We denote by $\AJ_X$ the Abel-Jacobi map of $X$, cf. \eqref{ajm};
by~$\gamma_n$, the Hermite constant, cf. Definition~\ref{53}; and by
$\stsys_1(g)$, the stable 1-systole of a metric $g$, cf. formula~\eqref{41b}.

\begin{theorem}[M. Gromov]
\label{11a} Let $X^n$ be a compact manifold of equal dimension and first Betti
number: $\dim(X)=b_1(X)=n$.  Then every metric $g$ on $X$ satisfies the following
optimal inequality:
\begin{equation}
\label{10} \deg(\AJ_X) \stsys_1(g)^n \leq \left( \gamma_n \right)^{\frac{n}{2}}
\vol_n(g).
\end{equation}
The boundary case of equality in inequality~\eqref{10} is attained by flat tori
whose group of deck transformations is a critical lattice in $\R^n$.
\end{theorem}
Note that the inequality is nonvacuous in the orientable case, only if the real
cuplength of $X$ is $n$, i.e., the Abel-Jacobi map $\AJ_X$ is of nonzero
algebraic degree.

The method of proof is to construct a suitable map $A$ from~$X$ to its Jacobi torus,
\begin{equation}
\label{75} A: X \to J_1(X).
\end{equation}
Here the torus is endowed with the flat Riemannian metric of least volume which
dominates the stable norm.  The map $A$ has the following properties:
\begin{enumerate}
\item[(a)]
$A$ belongs to the homotopy class of the Abel-Jacobi map $\AJ_X$;
\item[(b)]
the derivative of $A$ is volume-decreasing at every point.
\end{enumerate}
The existence of such a map $A$ follows from the techniques of \cite{BI94, BI95}, as
explained in \cite[pp.~259-260]{Gr3}. Inequality \eqref{10} now follows from the
definition of the Hermite constant.  This approach was systematized and generalized
in \cite{IK, BCIK2}, as explained in Section~\ref{ik}.

\section{Inequalities combining dimension and codimension 1}
\label{bank}

The material treated here is in \cite{BanK, BaKa}.  An optimal stable systolic
inequality, for $n$-manifolds $X$ with $b_1(X,\R)=1$, is due to J.~Hebda \cite{He}:
\begin{equation}
\label{hebda} \stsys_1(g) \sys_{n-1}(g)\le\vol_n(g),
\end{equation}
with equality if and only if $(X,g)$ admits a Riemannian submersion with connected
minimal fibers onto a circle.

Here the systole, $\sys_k$, and the stable systole, $\stsys_k$, coincide in
codimension 1 (i.e., when $k=n-1$) for orientable manifolds. Thus
$\sys_{n-1}=\stsys_{n-1}$.  In general, the $k$-systole is defined by
minimizing the $k$-dimensional area of {\em integral\/} $k$-cycles which are
not nullhomologous.

\begin{proof}[Proof of \eqref{hebda}]
Let $\omega \in H^1(X,\Z)_\R$ be a primitive element in the integer lattice in
cohomology, and similarly $\alpha\in H_1(X,\Z)_\R$ in homology.  Let~$\eta\in\omega$
be the harmonic 1-form for the metric $g$.  Then there exists a map $f$ to the
circle such that
\[
f: X\to S^1=\R/\Z, \quad df=\eta.
\]
Using the Cauchy-Schwartz inequality, we obtain
\begin{equation}
\begin{aligned}
\|\omega\|^*_2 (vol_n(g))^{1/2} & = |\eta|_{2}(vol_n(g))^{1/2} \\ & \geq \int_X |df|
dvol_n .
\end{aligned}
\end{equation}
Using the coarea formula, cf. \cite[3.2.11]{Fe1}, \cite[p.~267]{Ch2}, we obtain
\begin{equation}
\label{82}
\begin{aligned}
\displaystyle \|\omega\|^*_2 (vol_n(g))^{1/2} & \geq \int_{S^1} vol_{n-1} \left(
f^{-1}_{\phantom{I}}(t) \right)dt \\ & \geq \| \PD (\omega) \|_\infty \\ & \geq
\sys_{n-1}(g),
\end{aligned}
\end{equation}
since the hypersurface $f^{-1}(t) \subset X$ is Poincar\'{e} dual to $\omega$ for
every regular value $t$ of $f$.  Let us normalize the metric to unit total volume to
simplify the calculations.  In the case $b_1(X)=1$, we have
\begin{equation}
\label{83} \| \alpha \|_2 \; \| \omega \|^*_2 = 1.
\end{equation}
Therefore
\[
\begin{aligned}
\stsys_1(g) \sys_{n-1}(g) & = \|\alpha \|_\infty \; \sys_{n-1}(g) \\ & \leq \|\alpha
\|_\infty \; \| \omega \|^*_2 \\ & \leq \|\alpha \|_2 \; \| \omega \|^*_2 \\ & =1,
\end{aligned}
\]
proving Hebda's inequality~\eqref{hebda}.
\end{proof}

An optimal inequality, involving the conformal 1-systole \eqref{57}, is proved
in~\cite{BaKa}, namely inequality~\eqref{(5.2)} below.  This inequality
generalizes simultaneously Loewner's inequality \eqref{(1.1)}, Hebda's
inequality~\eqref{hebda}, the inequality \cite[Corollary 2.3]{BanK}, as well as
certain results of G. Paternain \cite{Pa}.  Let~$X$ be a compact, oriented,
$n$-dimensional manifold with positive first Betti number $b_1(X)\ge 1$.  Then
every metric~$g$ on~$X$ satisfies
\begin{equation}
\label{(5.2)} \confsys_1(g) \sys_{n-1}(g)\leq \gamma'_{b_1(X)^{\phantom{i}}}
\vol_n(g) ^{\frac {n-1} {n}},
\end{equation}
where equality occurs if and only if the following three conditions are satisfied:
\begin{enumerate}
\item[(a)]
the stable norm in $H_1(X;\R)$ is Euclidean;
\item[(b)]
the deck transformations of the Jacobi torus~\eqref{26} form a dual-critical
lattice~\eqref{22'};
\item[(c)]
the Abel-Jacobi map~\eqref{ajm} is a Riemannian submersion with connected minimal
fibers.
\end{enumerate}

The proof of inequality~\eqref{(5.2)} is similar to the calculation \eqref{82}.  To
conclude the argument, identity~\eqref{83} is replaced by the inequality
\begin{equation}
\label{85} \| \alpha \|_2 \; \| \omega \|^*_2 \leq \gamma'_{b_1(X)}
\end{equation}
for a pair of minimizing elements
\[
\alpha\in H_1(X,\Z)_\R \setminus \{ 0 \}, \quad \omega \in H^1(X,\Z)_\R \setminus \{
0 \}.
\]
Note that inequality~\eqref{85} results from the definition of the
Berg\'e-Martinet's constant $\gamma_b'$ combined with the fact that the $L^2$ norm
in cohomology is actually a Euclidean norm since harmonic forms form a vector space.

\section{Inequalities combining varieties of 1-systoles}
\label{ik}

A generalisation of Gromov's inequality~\eqref{10} was obtained in \cite{IK}, with a
number of further generalisations in \cite{BCIK2}. The following two theorems were
proved in \cite{IK} (the second one is an immediate corollary of the first one).
Denote by~$\fx \in H_{n-b}(\XX;A)$ the class of a typical fiber of the map
$\overline{\AJ}$ of~\eqref{72}, where $A=\Z$ or $\Z_2$ according as $\XX$ is
orientable or not.  Following Gromov \cite{Gr1}, denote by $\deg(\AJ_X)$ the infimum
of $(n-b)$-volumes of cycles representing the class~$\fx$.  Note that this quantity
is a topological invariant only if $n=b$.

\begin{theorem}[\cite{IK}]
\label{12b} Let $X$ be a closed Riemannian manifold of dimension $n=\dim(X)$.  Let
$b=b_1(X)$, and assume $n\ge b\geq 1$.  Then every metric $g$ on $X$ satisfies the
inequality
\begin{equation}
\label{eq12} {\deg(\AJ_X) \stsys_1(g)^b} \leq (\gamma_b) ^{\frac{b}{2}} {\vol_n
(g)}.
\end{equation}
\end{theorem}

\begin{theorem}[\cite{IK}]
\label{13} Let $X$ be a closed orientable manifold.  Let~$b=b_1(X)$.  Assume that
$\dim(X)=b+1$, and $\fx \not=0$.  Then every metric~$g$ on~$X$ satisfies the
following optimal inequality:
\begin{equation}
\label{11f} \stsys_1(g)^b \pisys_1(g) \leq (\gamma_b)^{\frac{b}{2}} \vol_{b+1}(g).
\end{equation}
\end{theorem}

The main tool in the proof of Theorems~\ref{12b} and \ref{13} is a generalisation of
the construction of an area-decreasing map \eqref{75} in a situation where the
dimension is bigger than the Betti number.  A suitable version of the coarea formula
is then applied to this map to conclude the proof.

The paper \cite{BCIK2} contains a number of generalisations of such stable
inequalities as well as the study of the boundary case of equality, with conclusions
along the closing lines of Section~\ref {bank}.  In particular, the Abel-Jacobi map
is a submersion.

Now let $X$ be a closed 3-manifold.  Let $\lambda$ be the extension of the
Casson-Walker invariant defined in~\cite{Les}.  The general definition is a bit
involved, but in our case the invariant can be expressed in terms of the
self-linking number of a typical fiber of the Abel-Jacobi map, cf. \eqref{102}.
We show in Section~\ref{thirteen} that the nonvanishing of the $\lambda$
invariant is a sufficient condition for the nonvanishing of the class $\fx$.

The argument runs roughly as follows.  Consider the map
$\overline{\AJ}: \XX \to \R^2$ of \eqref{72}.  Let $a,b\in \R^2$.  Let
$\bar F_a= \overline{\AJ}^{-1}(a)$ and $\bar F_b=
\overline{\AJ}^{-1}(b)$ be lifts of the corresponding fibers $F_a, F_b
\subset X$.  Choose a properly imbedded ray~$r_b\subset \R^2$ joining
$b$ to infinity while avoiding $a$ (as well as its translates), and
consider the complete surface \linebreak $S = \overline{\AJ}^{-1}(r_b)
\subset \XX$ with~$\partial S = \bar F_b$.  Furthermore, $S\cap
g. \bar F_a = \emptyset$ for all $g\in G$, where $G$ denotes the group
of deck transformations of the cover $p_G: \XX \to X$.  If $\bar
F_a$ is zero-homologous in $\XX$, choose a compact surface $M\subset
\XX$ with
\[
\partial M = \bar F_a.
\]
The linking number $\ell_X(F_a, F_b)$ in $X$ can therefore be computed
as the {\em algebraic\/} intersection
\[
\begin{aligned}
\ell_X(F_a, F_b) & = p_G(M) \cap F_b \\ & = \sum _{g\in G} g. M \cap
\bar F_b \\ & = \sum _{g\in G} \partial \left( g. M \cap S \right) \\
& = 0,
\end{aligned}
\]
where all sums and intersections are algebraic, and the curve $M \cap
S$ is compact by construction.  The details of a rigorous proof are
found in Section~\ref{marin}.

Therefore we obtain the following immediate corollary of Theorem~\ref{13}.

\begin{corollary}
\label{93b} 
Suppose $X$ is a three-manifold such that $b_1(X)=2$, with
nonvanishing $\lambda$-invariant.  Then every metric $g$ on $X$
satisfies the inequality
\begin{equation}
\label{93} \stsys_1(g)^2 \pisys_1(g) \leq \gamma_2 \vol(g).
\end{equation}
In the boundary case of equality, the manifold must be one of the bundles of
\linebreak Example~\ref{105}.
\end{corollary}
Analyzing the condition $\fx \not= 0$ leads to interesting questions about
3-mani\-folds, see Section~\ref{thirteen}.

\begin{remark}[Relation to LS]
In dimension 3, there is a connection between systoles and the
Lusternik-Schnirelmann category, denoted ${\rm cat}_{\rm LS}$.  A new
invariant ${\rm cat}_{\sys}$, called {\em systolic category}, is
defined in terms of the existence of systolic inequalities of a
suitable type~\cite{KR}.  It turns out that the two categories
coincide in dimension three: ${\rm cat}_{\rm LS}= {\rm
cat}_{\sys}$. This and other phenomena relating the two categories are
discussed in \cite{KR, KR2}.
\end{remark}

\begin{remark}[Codimension 2]
In the case $n-b=2$, inequality~\eqref{eq12} is related to the filling
area conjecture.  For simplicity, let $X=\R P^2 \times \T^2$.  It was
proved in \cite{IK} that every metric $g$ on $\T^2\times \R P^2$
satisfies a ``Pu times Loewner'' inequality
\begin{equation}
\label{ptl} \pisys_1(g)^2 \stsys_1(g)^2\leq \frac{\tilde \sigma_2 }{4}\gamma_2
\vol_4(g),
\end{equation}
where $\tilde \sigma_2$ is defined similarly to \eqref{42} by allowing
all possible topology of fillings.
\end{remark}

All the inequalities discussed so far depend on the rather restrictive
hypothesis $\dim(X)\geq b_1(X)$.  In the case $\dim(X) < b_1(X)$ it is harder
to obtain optimal inequalities, and typically the extremal metrics are not
smooth, see \cite{CK}. Recently, progress was achieved toward removing the
topological assumption in Loewner's inequality \eqref{(1.1)} (similarly to the
generalisation~\eqref{pu3} of Pu's inequality), which involves the study of
surfaces of genus $s\geq 2$.  Namely, all hyperelliptic surfaces satisfy
Loewner's inequality
\begin{equation}
\label{114} \pisys_1(X)^2 \leq \gamma_2 \area(X) \mbox{\ \ if $X$ hyperelliptic},
\end{equation}
and in particular \cite{FK} all metrics in genus 2
\cite{KS1}. Moreover, all surfaces of genus at least 20 satisfy
Loewner's inequality, as well \cite{KS2}.

\section{Fiber
class, linking number, and the generalized Casson invariant
\boldmath$\lambda$}
\label{thirteen}

\newtheorem{remars}[equation]{Remarks}
\newenvironment{remarks}{\begin{remars} \rm }{\end{remars}}

\newcommand{\G}{\`}
\newcommand{\CC}{{\mathbb C}}
\newcommand{\CX}{{\cal X}}
\newcommand{\CL}{{\cal L}}
\newcommand{\ZZ}{{\mathbb Z}}
\newcommand{\RR}{{\mathbb R}}
\newcommand{\QQ}{{\mathbb Q}}
\newcommand{\bp}{\noindent {\em Proof. }}
\newcommand{\eop}{\nopagebreak
            \hspace*{\fill}$\diamond$
            \medskip}
\newcommand{\topajast}{A_{\ast}}
\newcommand{\topaj}{A}
\newcommand{\topajd}{A_2}

The results of this section tie in with Corollary \ref{93b}.  Namely,
Proposition~\ref{propvoul} provides a necessary condition so that the invariant
$\lambda$ does not vanish.  The condition is the nonvanishing of the fiber class in
the maximal free abelian cover. Thus the nonvanishing of $\lambda$ is a sufficient
condition for the inequality~\eqref{93} to be satisfied.

From now on, unless otherwise mentioned, all the manifolds considered are smooth,
compact, and oriented.  Boundaries are oriented according to the outward normal
first convention.

In any (oriented) $3$-manifold $X$, the {\em linking number\/} of two disjoint
closed curves $\gamma$ and $d$ with null homology classes in $H_1(X;\QQ)$ is defined
as follows. Let $N(\gamma)$ be a tubular neighborhood of $\gamma$ disjoint from $d$.
There exists a curve $C(k\gamma)$ in $N(\gamma)$ that is homologous to $k\gamma$ in
$N(\gamma)$ for some integer $k>0$, and that bounds a surface $\Sigma(C(k\gamma))$
transverse to $d$. Then the {\em linking number\/} of $\gamma$ and $d$ in $X$ is
defined as
$$\ell_X(\gamma,d)=\frac{ \langle \Sigma(C(k\gamma)),d \rangle
_X}{k}$$ where $\langle \;,\; \rangle_X$ denotes the algebraic intersection number
in $X$.

\begin{proposition}
\label{propvoul} Let $X$ be a connected closed three-manifold whose first Betti
number is two.  Let $S_1$ and~$S_2$ be two surfaces imbedded in~$X$ whose homology
classes generate $H_2(X;\ZZ)$. The intersection of~$S_1$ and~$S_2$ is denoted by
$\gamma$. Let $\gamma^{\prime}$ be a parallel of $\gamma$ in $S_1$.  Let
$\overline{X}$ be the maximal free abelian cover of $X$, and let $\overline{\gamma}$
be a lift of~$\gamma$ in $\overline{X}$.  If the homology class of
$\overline{\gamma}$ vanishes in $H_1(\overline{X};\QQ)$, then the linking number
\begin{equation}
\label{lescop} \ell_X(\gamma,\gamma^{\prime})
\end{equation}
of $\gamma$ and $\gamma^{\prime}$ in $X$ vanishes, as well.
\end{proposition}

This proposition will be proved in Section~\ref{marin} as a particular case of a
more general result, Proposition~\ref{propal}, as suggested by Alexis Marin.

\begin{remark}
The two parallels of $\gamma$ in $S_1$ are homotopic to each other in the complement
of $\gamma$ in $X$ and they are also homotopic to the two parallels of $\gamma$ in
$S_2$. We could equivalently define $\gamma^{\prime}$ as the parallel of $\gamma$
with respect to the framing of $\gamma$ induced by $S_1$ or $S_2$.
\end{remark}

\begin{remark}
By Poincar\'e duality, under the assumptions of the proposition, there exists a
basis $(a_1, a_2)$ of $H_1(X;\QQ)$ such that the homology class $[d]$ of any closed
curve $d$ of $X$ reads
$$
[d]=\langle d,S_1\rangle_X^{\phantom{I}} a_1 + \langle
 d,S_2\rangle_X^{\phantom{I}} a_2,
$$
cf. \eqref{4.3}.  In particular, we have $[\gamma]=0$ in $H_1(X;\QQ)$.
Therefore, the number~$\ell_X(\gamma,\gamma^{\prime})$ is well-defined.  In
fact, if $\lambda$ denotes the extension of the Casson invariant defined in
\cite{Les}, and if $|\mbox{Torsion}(H_1(X))|$ is the cardinality of the torsion
part of $H_1(X;\ZZ)$, then
\begin{equation}
\label{102} \ell_X(\gamma,\gamma^{\prime}) =
-\frac{\lambda(X)}{|\mbox{Torsion}(H_1(X))|},
\end{equation}
cf. \eqref{103}.  Thus the linking number does not depend on the choice of the
transverse surfaces~$S_1$ and $S_2$ whose homology classes generate
$H_2(X;\ZZ)$. An easy direct proof of this fact is given in
\cite[pp.~93-94]{Les}.
\end{remark}

\begin{definition}
Let $X$ be a connected topological space, and let
\[
g: \pi_1(X) \to G
\]
be a homomorphism.  The {\em connected cover
\begin{equation}
\label{cc} \hat{X}(\mbox{Ker}(g))
\end{equation}
of $X$ associated to $g$\/} is the connected cover of $X$ where a loop of $X$ lifts
homeomorphically if and only if its homotopy class is in the kernel of~$g$.
\end{definition}
In particular, we have $\pi_1(\hat{X} (\mbox{Ker}(g) )) = \mbox {Ker} (g)$. The
covering group is the image of $g$.  For example, the maximal free abelian cover
$\XX$ of $X$ is associated with the natural composition
$$\pi_1(X) \longrightarrow H_1(X) \longrightarrow
H_1(X)/\mbox{Torsion}(H_1(X)).$$

Since $[\gamma]=0$ in $H_1(X;\QQ)$, the covering map $\rho$ from $\overline{X}$ to
$X$ maps~$\overline{\gamma}$ to $\gamma$, homeomorphically. In particular,
$\overline{\gamma}$ is a closed curve.

\section{Fiber class in NIL geometry}
\label{eleven}

Proposition~\ref{propvoul} may be illustrated by the following example
of the non-trivial circle bundles $N_e$ over the torus. For them, we
shall directly see that
\begin{itemize}
\item
  the roles of $\gamma$ and $\gamma^{\prime}$ can be played by two disjoint
fibers $F$ and $F^{\prime}$,
\item
  $\ell_{N_e}(F,F^{\prime})$ does not vanish, and
\item
  $H_1(\overline{N_e};\ZZ)=\ZZ \left[ \overline{F} \right]$.
\end{itemize}

\begin{example}
\label{105} Let $p_e:N_e \longrightarrow S^1 \times S^1$ be the circle bundle over
the torus with Euler number $e \neq 0$.
Denote by $x=S^1 \times \{1\}$ and $y=\{1\} \times S^1$ the two $S^1$-factors of the
base $S^1 \times S^1=x \times y$ of $N_e$.

Let~$D \subset x \times y$ be a disk disjoint from the union $x \cup y$.  The bundle
is trivial both over $D$ and over its complement. Thus we have
\[
N_e = S^1 \times s \left( \overline{(x \times y) \setminus D} \right) \cup_{S^1
\times \partial D} S^1 \times D ,
\]
where $s: \overline{(x \times y) \setminus D} \longrightarrow N_e$ is a section, and
the gluing map reads
\begin{equation}
\label{glue} (u,s(v)) \in S^1 \times (\partial D \cong S^1) \mapsto (uv^{\pm e},v)
\end{equation}
(up to signs).  Let $F=p_e^{-1}(x \cap y)$ be the fiber.  Then
$$
H_1(N_e)= \ZZ[s(x)] \oplus \ZZ[s(y)] \oplus \frac{\ZZ}{|e|\ZZ} [F].
$$
The surfaces $S_1$ and $S_2$ of Proposition~\ref{propvoul} may be chosen to be
$S_1=p_e^{-1}(x)$ and $S_2=p_e^{-1}(y)$, since these surfaces are dual to the basis
\[
\left( \left[ s(y) \right], \left[ s(x) \right] \right)
\]
of $H_1(N_e)/\mbox{Torsion}$ with respect to the algebraic intersection.  Then the
intersection $S_1 \cap S_2$ is precisely the fiber $F$, and its parallel induced by
the surfaces is another fiber $F^{\prime}$.

In this case, if the loop $s(\partial D)$ bounds a section of the bundle over
$s\left((x \!\times\! y) \setminus D\right)$, then the loop $(s(\partial D) \pm e
F)$ bounds a section of the bundle over $D$. This allows us to see that $|e|F$
bounds a surface that is pierced once by $F^{\prime}$, so that the linking number
satisfies
\begin{equation}
\label{103} \ell_{N_e}(F,F^{\prime})= \pm\; \frac{1}{e},
\end{equation}
cf. \eqref{102}.  Now, let us study the maximal free abelian covering
\begin{equation}
\label{104} \rho: \overline{N_e} \longrightarrow N_e
\end{equation}
associated with the map
\[
\pi_1(N_e) \longrightarrow (H_1(N_e)/\mbox{Torsion})=\ZZ [s(x)] \oplus \ZZ [s(y)]=G.
\]
The free abelian group $G$ acts on $\overline{N_e}$ freely as the covering
group. Therefore, $H_1(\overline{N_e})$ becomes a $\ZZ[G]$-module. The ring
$\ZZ[G]$ is denoted by $R$ and is identified with $\ZZ[t_x^{\pm 1}, t_y^{\pm
1}]$ by mapping $[s(x)]$ and $[s(y)]$ to~$t_x$ and $t_y$, respectively.
Consider  the standard product covering
$$
\begin{array}{llll}\rho_T:& \RR^2 &\longrightarrow &(S^1)^2\\
& (u,v) &\mapsto & (\exp(2i\pi u), \exp(2i\pi v)).\end{array}
$$
Then $\overline{N_e}$ is obtained from the trivialisation of $\rho$ of \eqref{104}
over the complement of $D$,
$$
\rho^{-1} \left( p_e^{-1} \left( (S^1)^2 \setminus D \right) \right) =F \times
\left( \RR^2 \setminus \rho_T^{-1}(D)\right),
$$
by filling in the $\ZZ^2$ holes (with boundaries the lifts of $F \times \partial D$)
by $\ZZ^2$ copies of $S^1 \times D$.  Here we use the same gluing map \eqref{glue}
as before,~$\ZZ^2$ times, equivariantly.  Then
\[
\begin{aligned}
H_1(\overline{N_e}) & =\frac{R [\overline{F}] \oplus R [\overline {\partial
D}]}{(t_x -1)R[\overline{F}] + (t_y -1) R [\overline {F}] + R \left(
[\overline{\partial D}] \pm e[ \overline {F} ] \right)} \\ & =\frac{R}{(t_x -1)R +
(t_y -1)R}[\overline{F}] \\ & =\ZZ \left[ \overline{F} \right].
\end{aligned}
\]
\end{example}

\section{A nonvanishing result (following A. Marin)}
\label{marin}

As it will be shown in Remark~\ref{rkvoulimpal} below, Proposition~\ref{propvoul} is
a particular case of Proposition~\ref{propal} that has been suggested and proved by
Alexis Marin.

\begin{proposition}
\label{propal} Let $X$ be a connected closed three-manifold.  Let~$S_1$ and $S_2$ be
two surfaces imbedded in $X$ that intersect transversely in $X$ along a curve
$\gamma$.
Let $\gamma^{\prime}$ be a parallel of $\gamma$ in $S_1$.  Let $\hat{X}_2$ be the
connected cover of $X$ associated with the composition
\begin{equation}
\label{111}
\begin{array}{llll} f_{2 \ast}:\pi_1(X) \longrightarrow & H_1(X)
&\longrightarrow &\ZZ\\ &[d] & \mapsto & \langle d,S_2\rangle_X
\end{array}
\end{equation}
where the first map is the Hurewicz homomorphism, cf. \eqref{cc}.  Let
$\hat{\gamma}$ be a lift of~$\gamma$ in $\hat{X}_2$.  If the homology class of
$\hat{\gamma}$ vanishes in $H_1(\hat{X}_2;\QQ)$, then the linking number
$\ell_X(\gamma,\gamma^{\prime})$ of $\gamma$ and $\gamma^{\prime}$ in $X$ is
well-defined and vanishes.
\end{proposition}

\begin{remark}
\label{rkvoulimpal} Let $\hat{X}$ be the connected cover of $X$ associated with the
composition
\begin{equation}
\label{112}
\begin{array}{llll}
\topajast:\pi_1(X) \longrightarrow & H_1(X) &\longrightarrow &\ZZ \oplus \ZZ\\ &[d]
& \mapsto & \left( \langle d,S_1\rangle_X,\langle d,S_2\rangle_X \right)
\end{array}
\end{equation}
where the first map is the Hurewicz homomorphism.  Then $\hat{X}$ is a cover of
$\hat{X}_2$.  Suppose the homology class of a lift of~$\gamma$ in $\hat{X}$ vanishes
in~$H_1(\hat{X};\QQ)$.  Then its image in $\hat{X}_2$ under the covering map, that
is a lift of~$\gamma$ in~$\hat{X}_2$, vanishes, as well.  Therefore, assuming that
Proposition~\ref{propal} is true, replacing $\hat{X}_2$ by $\hat{X}$ in its
statement yields another true proposition. In particular, Proposition~\ref{propal}
implies Proposition~\ref{propvoul} because in this case the above $\hat{X}$ is
precisely the maximal free abelian cover $\XX$.
\end{remark}

\begin{remark}
In Proposition~\ref{propal} that applies to $3$-manifolds with arbitrarily large
Betti numbers, the covering group is the subgroup $f_{2 \ast}(\pi_1(X))$ of $\ZZ$
that is isomorphic to $\{0\}$ or $\ZZ$.
\end{remark}

\begin{proof}[Proof of Proposition~\ref{propal}] Let $\rho_2:\hat{X}_2
\longrightarrow X$ be the covering map. Assume that the homology class of
$\hat{\gamma}$ vanishes in $H_1(\hat{X}_2;\QQ)$. Note that this implies that
$\gamma$ is rationally null-homologous in $X$ and that
$\ell_X(\gamma,\gamma^{\prime})$ is well-defined.

There exists a surface $\Sigma$ in $\hat{X}_2$ whose boundary $\partial \Sigma$ lies
in a small tubular neighborhood $N(\hat{\gamma})$ of $\hat{\gamma}$ such that
\begin{itemize}
\item $\partial \Sigma$ is homologous to $k[ \gamma]$ in $N(\hat{\gamma})$, for some $k>0$,
\item $N(\hat{\gamma})$ does not meet  $\rho_2^{-1}(\gamma^{\prime})$,
\item the restriction of $\rho_2$ to $N(\hat{\gamma})$ is injective.
\end{itemize}

Then
$$
k\ell_X(\gamma, \gamma^{\prime})=\langle \rho_2(\Sigma),\gamma^{\prime}\rangle_X
=\langle \Sigma, \rho_2^{-1}(\gamma^{\prime}) \rangle_{\hat{X}_2}\;
$$
and $k\ell_X(\gamma, \gamma^{\prime})$ is the sum over the lifts of
$\gamma^{\prime}$ of their algebraic intersections with $\Sigma$.

We now apply the following Lemma~\ref{lem114} to complete the proof of
Proposition~\ref{propal}.
\end{proof}

\begin{lemma}
\label{lem114} For any lift $\hat{\gamma}^{\prime}$ of $\gamma^{\prime}$, the
algebraic intersection of $\Sigma$ and~$\hat{\gamma}^{\prime}$ vanishes.
\end{lemma}
\begin{proof}[Proof of Lemma~\ref{lem114}]
The proof consists in seeing $\langle \Sigma, \hat{\gamma}^{\prime}
\rangle_{\hat{X}_2}$ as the algebraic boundary of the compact oriented intersection
of $\Sigma$ with a noncompact surface $\hat{S}_1^*$ that is bounded by
$\hat{\gamma}^{\prime}$ and that does not meet~$N(\hat{\gamma})$.

Let us first construct $\hat{S}_1^*$.

Consider a tubular neighborhood $S_2 \times [-1,1]$ of $(S_2=S_2 \times {0})$
such that\linebreak $ S_1 \cap (S_2 \times [-1,1])$ reads $\gamma \times
[-1,1]$ in $S_2 \times [-1,1]$ (up to orientations).

Let  $f_2: X \to S^1$ be the map to the unit circle in $\C$, defined as follows on
the subset~$S_2 \times [-1,1] \subset X$:
\[
\begin{array}{llll}f_2:&S_2 \times [-1,1] &\longrightarrow & S^1\\
&(x,t) & \mapsto & \exp(i\pi t), \end{array}
\]
while $f_2$ maps the complement of $S_2 \times [-1,1]$ in $X$ to the point $(-1)$.

Note that the map~$f_2$ induces the morphism~$f_{2 \ast}$ of~\eqref{111}. Thus the
composition~$f_2 \circ \rho_2$ induces the trivial map from $\pi_1(\hat{X}_2)$ to
$\pi_1(S^1)$. Therefore, it factors through the covering
$$\begin{array}{llll}\mbox{Exp}:& \RR &\longrightarrow &S^1\\ & u
&\mapsto & \exp(2i\pi u)\end{array}$$ yielding a map $\hat{f}_2 :\hat{X}_2
\longrightarrow \RR$ such that $\mbox{Exp} \circ \hat{f}_2=f_2 \circ \rho_2$.

Without loss assume that $$\gamma^{\prime}=S_1 \cap f_2^{-1}(i)$$ and that
$$N(\hat{\gamma}) \subset \hat{f}_2^{-1}
\left[-\tfrac{1}{8},\tfrac{1}{8}\right].$$

Now, any lift $\hat{\gamma}^{\prime}$ of $\gamma^{\prime}$ reads
\begin{equation}
\rho_2^{-1} (S_1) \cap \hat{f}_2^{-1} (\{n+1/4\})
\end{equation}
for some $n \in \ZZ$, and splits~$\rho_2^{-1} (S_1)$ as the union of
$$
\hat {S}_1^+= \rho_2^{-1} (S_1) \cap \hat{f}_2^{-1} \left(n+ \tfrac{1}{4} +
[0,+\infty[ \right)
$$
and
$$
\hat{S}_1^- =\rho_2^{-1} (S_1) \cap \hat{f}_2^{-1} \left( n+ \tfrac{1}{4} + \;
]\!-\infty,0] \right).
$$
One of the non-compact subsurfaces $\hat{S}_1^+$ and $\hat{S}_1^-$ does not
meet~$N(\hat{\gamma})$. This will be our surface $\hat{S}_1^*$.

We may assume that the surface $\Sigma$ is transverse to $\hat{S}_1^*$. Then these
two surfaces intersect along an oriented curve $\Sigma \cap \hat{S}_1^*$ that is
oriented so that the triple
\begin{enumerate}
\item[(a)]
tangent vector $\vec{t}(\Sigma \cap \hat{S}_1^*)$ to the intersection curve,
\item[(b)]
positive normal vector $\vec{n}(\Sigma)$ to $\Sigma$,
\item[(c)]
positive normal vector $\vec{n}(\hat{S}_1^*)$ to~$\hat{S}_1^*$
\end{enumerate}
is direct.  The curve $\Sigma \cap \hat{S}_1^*$ is a collection of arcs and closed
curves properly imbedded in $\hat{S}_1^*$ since $\hat{S}_1^*$ does not
meet~$\partial \Sigma$.  Since the oriented boundary of $\Sigma \cap \hat{S}_1^*$
(that obviously vanishes in $H_0(X)$) represents the algebraic intersection of
$(\partial \hat{S}_1^*=\pm \hat{\gamma}^{\prime})$ with $\Sigma$, the lemma is
proved.

Let us be slightly more explicit about this last argument.  A point~$x$ of the
oriented boundary of $\Sigma \cap \hat{S}_1^*$ gets the sign $\varepsilon=\pm 1$
such that the tangent vector $\vec{t}_x(\Sigma \cap \hat{S}_1^*)$ of $\Sigma \cap
\hat{S}_1^*$ at $x$ is oriented as $\varepsilon \vec{n}_x(\Sigma \cap \hat{S}_1^*)$
where $\vec{n}_x(\Sigma \cap \hat{S}_1^*)$ is the outward normal vector to the curve
at $x$.  (In other words, final points of arcs get a positive sign while initial
points get a minus sign.)  On the other hand, $\vec{n}_x(\Sigma \cap \hat{S}_1^*)$
may be identified to the outward normal of $\hat{S}_1^*$ and the triple
\[
\left( \vec{n}_x \left( \Sigma \cap \hat{S}_1^* \right), \vec{t}_x \left( \partial
\hat {S} _1 ^* \right), \vec{n}_x \left( \hat{S}_1^* \right) \right)
\]
is direct. Furthermore, $\vec{n}_x(\Sigma)$ may be identified with
$\mbox{sign}(x)\vec{t}_x(\partial \hat{S}_1^*)$ where $\mbox{sign}(x)$ is the sign
of the intersection of $\partial \hat{S}_1^*$ and $\Sigma$ at $x$.  Thus, the triple
\[
\left( \varepsilon \vec{n}_x \left( \Sigma \cap \hat{S}_1^* \right), \; \mbox{sign}
(x) \vec{t}_x \left( \partial \hat{S}_1^* \right), \; \vec{n}_x \left( \hat{S}_1^*
\right) \right)
\]
is direct, too.  This shows that $\varepsilon=\mbox{sign}(x)$.  Then the sum of the
$\varepsilon$ vanishes because there are as many initial points of arcs in $\Sigma
\cap \hat{S}_1^*$ as there are final points, and that makes the algebraic
intersection of $\hat{\gamma}^{\prime}$ and $\Sigma$ that is the sum of the
$\mbox{sign}(x)$ vanish, too.  This concludes the proof of Lemma~\ref{lem114} and
the proof of the proposition.
\end{proof}

Note that Lemma~\ref{lem114} yields the following corollary.

\begin{corollary} Let $X$ be a $3$-manifold with Betti
number $2.$ Suppose the lift $L$ of a typical fiber of the Abel-Jacobi map, to
an infinite cyclic cover $C$ of $X$ is rationally zero-homologous. Equip $L$
with the framing induced by  the Abel-Jacobi map.
 Then
the self-linking of $L \subset C$ is zero.
\end{corollary}

This corollary is obviously not true when $C$ is replaced by the compact
manifold~$X$ itself, or by any {\em finite\/} cyclic cover of $X$.

\section*{Acknowledgments}
We are grateful to A. Marin for suggesting and proving Proposition~\ref{propal}.


\bibliographystyle{amsplain}

\end{document}